\documentclass[11pt]{article}
\usepackage[dvips]{graphicx}
\usepackage{a4wide}
\usepackage{amsmath}
\usepackage{amsthm}
\usepackage{amsfonts}
\usepackage{bm}
\usepackage{cases}
\usepackage{mathrsfs}
\usepackage{authblk}

{\theoremstyle{plain}%
  \newtheorem{thm}{\bf Theorem}[section]%
  \newtheorem{lem}[thm]{\bf Lemma}%
  \newtheorem{prop}[thm]{\bf Proposition}%
  \newtheorem{claim}[thm]{\bf Claim}%
}
{\theoremstyle{plain}
  \newtheorem{rem}[thm]{\bf Remark}
}

\newcommand{\dC}{{\mathbb{C}}}%
\newcommand{\dN}{{\mathbb{N}}}%
\newcommand{\dZ}{{\mathbb{Z}}}%
\newcommand{\dR}{{\mathbb{R}}}%
\newcommand{\dQ}{{\mathbb{Q}}}%

\title{Periodicity of Grover walks on distance-regular graphs}

\author[1]{Yusuke YOSHIE\footnote{E-mail: yusuke.yoshie.r1@dc.tohoku.ac.jp}}
\affil[1]{Graduate School of Information Sciences, Tohoku University, \newline Aoba, Sendai 980-8579, Japan}
\date{}
\begin {document}
\maketitle{}
\begin{abstract}
Characterizations graphs of some classes to induce periodic Grover walks have been studied for recent years. In particular, for the strongly regular graphs, it has been known that there are only three kinds of such graphs.  Here, we focus on the periodicity of the Grover walks on distance-regular graphs. The distance-regular graph can be regarded as  a kind of generalization of the strongly regular graphs and the typical graph with an equitable partition. In this paper, we find some classes of such  distance-regular graphs and obtain some useful necessary conditions to induce periodic Grover walks on the general distance-regular graphs. Also, we apply this necessary condition to give another proof for the strong regular graphs.
\end{abstract}
{\bf keywords: } Grover walk; Periodicity; Distance-regular graph;

\section{Introduction}
Quantum walks were introduced as quantizations of random walks \cite{Gudder} and these have been investigated for the last decade \cite{Konno14}, \cite{Wang}, \cite{Portugal}. Quantum walks can be defined as a model in which the walker moves around an underlying graphs like random walks but the behavior of quantum walks are definitely different from that of random walks. For random walks, the probability distribution at each time can be given by a vector on a space induced by the underlying graph and time evolution can be given by the transition matrix. For quantum walks, the dynamics of the walker can be interpreted as a motion of plane wave on the graph \cite{HKSS13}. So the walker has amplitude  expressed by a complex value.  Then the state of the walker at each time can be expressed by a vector called {\it quantum state}  and the time evolution can be obtained by iterating a unitary operator induced by the underlying graph. However, in order to gain the probability, we take the square of the absolute value of the quantum state. By such a non-liner calculation, the distribution of probability cannot be expressed by that of past times. Due to these differences, individual phenomena of quantum walks occur, e.g, coexistence of the localization \cite{Inui04} and linearly spreading, perfect state transfer \cite{Kendon}, \cite{Stefanak} and so on.  On these results, such phenomena can be observed in a fixed graph. On the contrary, our aim in this paper is to consider a kind of inverse problem from the viewpoint of quantum walks, that is, we first classify the phenomenon of quantum walks and find a graph in which the phenomenon occurs. Here, what we choose as such a phenomenon is the {\it periodicity}. The periodicity means that there exists an integer $k$ such that the $k$-th power of the evolution operator becomes the identity operator. If it happens, then the quantum state at time $k$ returns to the initial quantum state and the behavior becomes periodic. Since the underlying graph induces the time evolution operator, the graph determines the quantum walk is periodic or not, and we are interested in characterizing such graphs. Since the quantum walk is periodic if and only if all the eigenvalues of its evolution operator are roots of unity, spectral analysis is one of useful tool to consider the periodicity. In this paper, we refer the periodicity of the Grover walk. Indeed, the spectrum of the time evolution operator of the Grover walk can be obtained by projecting that of the transition matrix of simple random walk of the underlying graph onto unit circle. So analysis of spectrum of simple random walk is one of methods to consider the periodicity of the Grover walk. Taking this method, characterizations of such graphs have been obtained in \cite{Higuchi13}, \cite{Konno13}, \cite{Yoshie}. In particular, we take key results in this paper from these references. 

\subsection{Related works}
Here, we introduce related works for the periodicity of the Grover walks. Recently, characterizations of graphs with an {\it equitable partition} \cite{RGodsil} to induce periodic Grover walks have been found. For example, the perfect characterizations of the generalized Bethe tree, which is a tree graph with an equitable partition, to induce periodic Grover walks is referred in \cite{YoshieB}. Such Bethe trees are 
only the subdivision of two graphs arranged in Figures \ref{fig bethe3}, \ref{fig bethe4}. 
Also, we give the results about strongly regular graphs. The strongly regular graph $\mathrm{SRG}(n, k, \lambda, \mu)$ is a special kind of regular graphs with the following properties: (i) the number of vertices is $n$, (ii) the valency is $k$, (iii) any two adjacent vertices have $\lambda$ co-neighbors, (iv) any two non-adjacent vertices have $\mu$ co-neighbors. It also has an equitable partition. Then the following Proposition is known. 
\begin{prop} (Higuchi, Konno, Sato, Segawa \cite{Higuchi13})\\
The Grover walks induces by strongly regular graphs $\mathrm{SRG}(n, k, \lambda, \mu)$ are periodic if and only if $(n, k, \lambda, \mu)=(5, 2, 0, 1), (2k, k, 0, k), (3\lambda ,2\lambda ,\lambda ,2\lambda)$, whose periods are $5, 4, 12$, respectively. \label{prop srg}
\end{prop}
Then the graphs listed in Proposition \ref{prop srg} are nothing but $C_{5}, K_{k, k}, K_{\lambda, \lambda, \lambda}$, respectively. In fact, the strongly regular graph is a distance-regular graph whose diameter is $2$. So the distance-regular graphs can be regarded as generalized strongly regular graphs. The distance-regular graphs also have an equitable partition and these are typical examples of such graphs. So we treat the periodicity of the Grover walk on a typical example of such graphs in this paper. 
\begin{figure}[h]
\hspace{2cm}
\begin{tabular}{ccc}
		\begin{minipage}[h]{0.25\linewidth}
		\centering
		\vspace{1.5cm}
		\scalebox{1.0}{
		\includegraphics[scale=1.0]{bethe-figfig.1}
		}
		\vspace{0.1cm}
		\hspace{15cm}
             \centering
            \vspace{0.5cm}
		\caption{}
		\label{fig bethe3}
		\end{minipage}
		\hspace{1cm}
		\begin{minipage}[h]{0.25\linewidth}
		\centering
		\scalebox{1.0}{
		\includegraphics[scale=1.0]{bethe-figfig.2}
		}
		\caption{}
		
		\label{fig bethe4}
		\end{minipage}	
\end{tabular}
\end{figure}

\subsection{Main results}
In this paper, we say that $G$ is {\it periodic} if it induces a periodic Grover walk. For a graph $G$, let $T=T(G)$ be the transition matrix for simple random walk on $G$ (detailed definition is referred in Section 2). First, we give a necessary condition for a periodic graph with rational eigenvalues. 
\begin{thm}
For a periodic graph $G$, if all the eigenvalues of $T(G)$ are rational, then it should hold
\[ \mathrm{diam}(G) < 5, \] \label{thm of diam}
\end{thm}
where $\mathrm{diam}(G)$ is the diameter of $G$. So we can control the structure of graph by the condition of periodicity. Next, we focus on two well-known distance regular graphs, Hamming graphs and Johnson graphs, which can be determined by two parameters. (These definitions are referred in Section 3). Then the followings are our results. 
\begin{thm}
The Hamming graph $H(d, q)$ is periodic if and only if 
\[(d, q)\in \{ (1, 2), (1, 3), (2, 2), (3, 3), (4, 2)\}. \] \label{thm of hamming}
\end{thm}
\begin{thm}
The Johnson graph $J(n, k)$ is periodic if and only if 
\[ (n, k)\in \{ (2, 1), (3, 1), (4, 2) \}. \] \label{thm of johnson}
\end{thm}
In addition, we treat periodicities of the Grover walks of general distance-regular graphs. For a distance-regular graph $G$ with diameter $d$, we define $\varphi(x)$ as
\begin{equation}
\varphi(x):=\mathrm{det}(x I_{D}-\tilde{Q})=\sum^{D}_{i=0}\rho_{i}x^{i}, \label{def of phi}
\end{equation} 
where $D=d+1$ and $I_{D}, \tilde{Q}$  are the $D \times D$ identity matrix,  the modified quotient matrix of $G$, respectively (detailed definition is found in Section 4). Then the following holds. 
\begin{thm}
For a general distance-regular graph $G$, let $\rho_{i}, D$ be defined as the previous arguments. If $G$ is periodic, then it should hold
\[ 2^{j}\rho_{D-j} \in \dZ \]
for every $j \in \{1, \cdots, D-1 \}$. \label{thm of rho}
\end{thm}

This paper is organized as follows: In section 2, we define the Grover walk on the graphs and consider its periodicity and we give the definition of distance-regular graphs and introduce the spectrum of these graphs. In section 3, we prove Theorems \ref{thm of hamming}, \ref{thm of johnson}. Next, we refer applications to the general distance-regular graphs and obtain some necessary conditions for the general distance-regular graphs to induce periodic Grover walks in section 4. Section 5 is devoted to summarize our results and make a discussion for our future works. 
\section{Preliminaries}
 \subsection{Periodicity of the Grover walk}
Here, we define the Grover walk on graphs and consider its periodicity. Let $G=(V(G), E(G))$ be a finite simple graph with the vertex set $V(G)$, and the edge set $E(G)$. For $uv \in E(G)$, the arc from $u$ to $v$ is denoted by $(u, v)$. The origin and terminus vertices of $e=(u, v)$ are denoted by $o(e), t(e)$, respectively and we express $e^{-1}$ as the inverse arc of $e$. We define $\mathcal{A}(G)=\{ (u, v), (v, u) | uv \in E \}$, which is the set of the symmetric arcs of $G$. Then the Grover transfer matrix of $G$ is defined by the following unitary matrix $U=U(G)$ indexed by $\mathcal{A}(G)$:
\begin{equation*}
U_{e, f}=
	\begin{cases}
	2/\mathrm{deg}(t(f)) & \text{if $t(f)=o(e), e\neq f^{-1}$,}\\
        2/\mathrm{deg}(t(f))-1 & \text{if $e=f^{-1}$,}\\
        0 & \text{otherwise.}
	\end{cases}
\end{equation*}  
Then it is given as the time evolution operator of the Grover walk on the graph $G$. Also, let $\varphi_{t} \in \ell^{2}(\mathcal{A}(G))$ be a quantum state at time $t$, which is a vector representing amplitude of each arcs at time $t$. The finding probability on $e$ at time $t$ is defined by $|\varphi_{t}(e)|^{2}$, which is the square of absolute value of the entry corresponding to arc $e$ of $\varphi_{t}$.  Moreover, we can obtain $\varphi_{t}$ by $t$-th iteration of $U$ to the initial quantum state $\varphi_{0}$, that is, $\varphi_{t}=U^{t}\varphi_{0}$. Due to its unitarity, the norm is preserved. Since the Grover transfer matrix is determined by the underlying graph, we can say that the graph induces the Grover walk. So the property of the Grover walk depends on the underlying graph. In particular, we introduce the periodicity of the Grover walk and treat some special classes of graphs to induce periodic Grover walks. For $k \in \dN$, we say that a graph $G$ induces a $k$-periodic Grover walk if and only if $U^{k}=I_{|\mathcal{A}(G)|}$ and $U^{j} \ne I_{|\mathcal{A}(G)|}$ for every $j$ with $0<j<k$. So it immediately follows that in $k$-periodic graphs,  $\varphi_{k}$ returns to $\varphi_{0}$ for {\it any} initial quantum state $\varphi_{0}$. For a square matrix $X$, the set of the eigenvalues of $X$ is denoted by $\sigma(X)$. Here, we give an equivalent Proposition for $k$-periodic graphs, which plays an important role in this paper: 
\begin{prop}
A graph $G$ is a $k$-periodic graph if and only if $\lambda^{k}_{U}=1$ for every $\lambda_{U} \in \sigma(U)$, and there exists $\lambda_{U} \in \sigma(U)$ such that $\lambda^{j}_{U} \ne 1$ for every $j$ with $0 < j < k$.\label{prop of p}
\end{prop}
Therefore, we use a spectral method in order to analyze the periodicity of the Grover walk with the above Proposition. Let $T=T(G)$ be the transition matrix for simple random walk on $G$, that is, for $u, v \in V$, 
\begin{equation*}
(T)_{u, v}
=
\begin{cases}
	1/\mathrm{deg}(u) & \text{if $u \sim v$,}\\
	0 & \text{otherwise.}
	\end{cases}
\end{equation*}
Since all the summations over each rows of $T$ is $1$, the $|V|$-dimensional all-one vector is always the eigenvector for eigenvalue $1$ and it gives the maximum eigenvalue. So the absolute values of all the eigenvalues of $T$ should be at most $1$. Results for the eigenvalues of $U$ have been found in, e.g. \cite{Emms06}, \cite{HSegawa13} and the above $T$ is strongly related to the spectrum of the Grover walk.  
\begin{lem}
(Higuchi, Segawa \cite{HSegawa17})\\
The spectrum of the Grover transfer matrix $U$ is decomposed by
\[\sigma(U)= \{ e^{\pm i \arccos{ \left( \sigma(T) \right)} } \}  \cup {\{ 1 \}}^{b_{1}}  \cup  {\{ -1 \} }^{b_{1}+1-{\bf 1_{B}}},  \]
where $b_{1}$ is the first Betti number of $G$, that is, $|E(G)|-|V(G)|+1$, and ${\bf 1_{B}}= 1$ if G is bipartite,  ${\bf 1_{B}}= 0$ otherwise. 
\end{lem}
Therefore, if a graph $G$ is periodic, then all the eigenvalues of $T$ should be the real parts of a root of unity. Then the following Zhukovskij transformation can be useful tool for our researches .  
\begin{lem}
(Higuchi, Konno, Sato, Segawa \cite{Higuchi13})\\
Let $f(\lambda)$ be a monic polynomial of degree $i$ for $\lambda \in \dR$. Then the solutions of $f(\lambda)=0$ are the real parts of  roots of unity if and only if for $z \in \dC$ with $|z|=1$, the polynomial  $(2z)^{i}f\left( (z+z^{-1})/2 \right)$ becomes a product of some cyclotomic polynomials. \label{lem of cyclo}
\end{lem}

\subsection{Definition of distance-regular graphs}
Here, we give the definition of the distance regular graphs. Let $G$ be a finite connected $k$-regular graph. The graph distance between $x, y \in V(G)$ is denoted by $d(x, y)$, which is the shortest path from $x$ to $y$. Also, we denote the set of neighbors of $x \in V(G)$ by $N(x)$ and the diameter of $G$, $d=\mathrm{diam}(G)$ is defined by $d=\mathrm{max}\{d(x, y) | x, y \in V(G) \}$. For a vertex $x \in V(G)$, we define $\Gamma_{j}(x):=\{ y \in V(G)\, |\, d(x, y)=j \}$, where $0 \le j \le d$. Then $G$ is a {\it distance-regular graph} if for any $x \in V(G)$ and $j$ with $0 \le j \le d$, the values $|N(y)\cap \Gamma_{j-1}(x)|, |N(y)\cap \Gamma_{j}(x)|, |N(y)\cap \Gamma_{j+1}(x)|$ are constant only depending on $j$ whenever $y \in \Gamma_{j}(x)$, that is, 
\begin{eqnarray}
|N(y)\cap \Gamma_{j-1}(x)|=&c_{j}, \\
|N(y)\cap \Gamma_{j}(x)|=&a_{j}, \\
|N(y)\cap \Gamma_{j+1}(x)|=&b_{j}, 
\end{eqnarray}
for any $x \in V(G)$ and $y \in \Gamma_{j}(x)$. Moreover,  $\{ \Gamma_{0}(x), \Gamma_{1}(x), \cdots, \Gamma_{d}(x) \}$ becomes an equitable partition of $G$. These positive parameters $a_{j}, b_{j}, c_{j}$ are called an {\it intersection array} \cite{GM} and it should hold that for $j \in \{0, 1, \cdots, d \}$, 
\begin{equation}
c_{j}+a_{j}+b_{j}=k, \label{c+a+b}
\end{equation}
where we define $c_{0}=b_{d}=0$. From the connectivity of $G$, it immediately follows that $b_{i}\ne 0, c_{j}\ne 0$ for $i \in \{0, 1, \cdots, d-1 \}$, $j \in \{1, 2, \cdots, d \}$.
Then we can induce the following $(d+1) \times (d+1)$ tri-diagonal matrix $Q$ called the {\it quotient matrix} of $G$, 
\begin{equation*}
Q=\left(
	\begin{array}{ccccccc}
	a_{0} & b_{0} & & &  & \\
	c_{1} & a_{1} & b_{1} & & & \\
	 & c_{2} & a_{2} & b_{2} & & \\
	 &  & \ddots & \ddots & \ddots & \\
	 & & & c_{d-1} & a_{d-1} & b_{d-1} \\
	 & & & & c_{d} & a_{d} \\
	 \end{array}
	\right).
\end{equation*}
Then it is well-known that $\sigma(A)=\sigma(Q)$, where $A$ is the adjacency matrix of $G$ \cite{GM}. Since $G$ is a $k$-regular graph, it follows
\[ \sigma(T)=\sigma \left( \frac{1}{k}Q \right). \]
Let us put $\tilde{Q}:=(1/k)Q$, which is a transition matrix on the projected path graph of $G$. So if a distance-regular graph is periodic, then the eigenvalues of its modified quotient matrix $\tilde{Q}$ should be the real parts of  a root of unity. 
\section{Periodicities of Grover walks on Hamming graph and Johnson graph}
\subsection{Proof of Theorem \ref{thm of diam}}
For a square matrix $X$, the number of distinct eigenvalues of $X$ is denoted by $n(X)$ in this paper. Then for any graph $G$ and its adjacency matrix $A=A(G)$, it is well-known that $\mathrm{diam}(G)<n(A)$ \cite{Spg}. Indeed, for its transition matrix $T=T(G)$, this relation similarly holds, that is, 
\begin{equation}
\mathrm{diam}(G)<n(T), \label{diam t}
\end{equation}
which can be proven by the same way. Using this fact, we can prove Theorem \ref{thm of diam}. 

\begin{proof}[Proof of Theorem \ref{thm of diam}]
For any $\lambda_{T} \in \sigma(T)$, it holds $\lambda_{T}=p/q$ for some $p, q \in \dN$ with $\mathrm{gcd}(p, q)=1$ and $|p/q|<1$ from the assumption. Put $h(x)=x-p/q$, which is the minimum polynomial of $\lambda_{T}$ on $\dQ$. Then we have
\begin{equation}
(2z)h\left(\frac{z+z^{-1}}{2} \right)=z^{2}-\frac{2p}{q}z+1. \label{zhuko of h}
\end{equation}
Since $G$ is periodic, $\lambda_{T}$ should be a real part of the root of unity, that is, for $z \in \dC$ with $|z|=1$, (\ref{zhuko of h}) should be represented by a product of cyclotomic polynomials by Lemma \ref{lem of cyclo}. So the candidates of $p/q$ are only $0, \pm1$, or $\pm 1/2$ since the above polynomial should be an integer polynomial. Indeed, one can easily show that (\ref{zhuko of h}) becomes a product of cyclotomic polynomials for these candidates. Under the assumption, it follows $\lambda_{T}=0, \pm 1$, or $\pm 1/2$ and  $n(T) \le 5$. Therefore, we have $\mathrm{diam}(G) < 5$ by (\ref{diam t}). 
\end{proof}
This Theorem can be applied to any graphs. In particular, it plays an important role for some classes of distance-regular graphs. 

\subsection{Hamming graphs}
For positive integers $d, q$, the Hamming graph $H(d, q)$ is defined as follows: Let $F$ be a finite set of $q$ elements. Then, the vertex set of $H(d, q)$ is $F^{d}$ and two vertices $x=(x_{1}, x_{2}, \cdots, x_{d}), y=(y_{1}, y_{2}, \cdots, y_{d}) \in F^{d}$ are adjacent if and only if the value $|\{ i \,| \,x_{i} \ne y_{i}, 1\le i \le d \}|=1$. Indeed, it is a $d(q-1)$-regular graph and $\mathrm{diam}(H(d, q))=d$. Let $A, T$ be the adjacency matrix and the transition matrix of $H(d, q)$, respectively. It is known that the distinct eigenvalues of its adjacency matrix are written by 
\[ \sigma(A)=\{ d(q-1)-qi \,| \,0 \le i \le d \} \]
\cite{GM}. So it immediately follows that
\[ \sigma(T)=\left\{ 1-\frac{qi}{d(q-1)} \,\Big{|} \,0 \le i \le d \right\}. \]
Using these facts, we prove Theorem \ref{thm of hamming}. 
\begin{proof}[Proof of Theorem \ref{thm of hamming}]
For the case $d=1$, the Hamming graph is nothing but the complete graph $K_{q}$. So in this case, $q$ should be only $2$ or $3$ \cite{Higuchi13}. Then the periodic Hamming graphs with $d=1$ are $K_{2}, K_{3}$. Therefore, we can suppose $d \ge 2$. Then all of the eigenvalues of $T$ are rational. So it follows $d < 5$ by Theorem \ref{thm of diam}. For the case $d=2$, the possibility of the spectrum of  $T$ is only 
\[ \sigma(T)=\left( 1-\frac{2q}{2(q-1)}, 1-\frac{q}{2(q-1)}, 1 \right)=(-1, 0, 1) \]
for $q=2$ since these rational eigenvalues should be $\pm 1/2$, or $\pm 1$, or $0$ and $q \in \dN$. Thus, the periodic Hamming graph is only $H(2, 2)\simeq C_{4}$ in this case. For the case $d=3$, the possibility is only 
\[ \sigma(T)=\left( 1-\frac{q}{(q-1)}, 1-\frac{2q}{3(q-1)}, 1-\frac{q}{3(q-1)}, 1 \right)=\left(-\frac{1}{2}, 0, \frac{1}{2}, 1\right) \]
for $q=3$. Thus, the graph is only $H(3, 3)$ in this case. For the case $d=4$ there is the only possibility such that 
\[ \sigma(T)=\left( 1-\frac{q}{(q-1)}, 1-\frac{3q}{4(q-1)}, 1-\frac{q}{2(q-1)}, 1-\frac{q}{4(q-1)}, 1 \right)=\left(-1, -\frac{1}{2}, 0, \frac{1}{2}, 1 \right) \]
for $q=2$. Thus, the graph is only $H(4, 2)$ in this case. Therefore, all the periodic Hamming graphs are only $H(1, 2), H(1, 3), H(2, 2), H(3, 3), H(4, 2)$.
\end{proof}

\subsection{Johnson graphs}
For two positive integers $n, k$ with $n \ge k$, the Johnson graph $J(n, k)$ is defined as follows: The  vertices of $J(n, k)$ are the $k$-element subsets of an $n$-element set. Also, two vertices $x, y \in V(J(n, k))$ are adjacent if and only if $|x\cap y|=k-1$. So this graph is a $k(n-k)$-regular graph and $\mathrm{diam}(J(n, k))=\mathrm{min}\{ k, n-k \}$. Let $A, T$ be the adjacency matrix and the transition matrix of $J(n, k)$, respectively. Also, it is known that the $(d+1)$ distinct eigenvalues of $A$  are written  by
\[ \sigma(A)=\{ (d-j)(n-d-j)-j \,| \,0 \le i \le d \}, \]
where $d=\mathrm{diam}(J(n, k))$ \cite{GM}. Therefore, we can obtain 
\[ \sigma(T)=\left\{ \frac{(d-j)(n-d-j)-j}{d(n-d)} \,\Big{|}\,0 \le j \le d \right\}. \]
\begin{proof}[Proof of Theorem \ref{thm of johnson}]
Similarly to the case of the Hamming graphs, all the eigenvalues of the transition matrix of $J(n, k)$ are rational. Thus, we have $0 < d < 5$ by Theorem \ref{thm of diam}.  For the case $d=1$, the possibilities of the spectrum of $T$ are only
\[ \sigma(T)=\left( -\frac{1}{n-1}, 1 \right)=\left(-1, 1\right), \left(-\frac{1}{2}, 1\right). \]
So we have $n=2$ or $3$ in this case. If $d=k$, that is,  $2k \le n$, then we can obtain $k=1$ for $n=2, 3$.  If $d=n-k$, that is, $2k \ge n$, then we can obtain $k=1, 2$ for $n=2, 3$, respectively. Indeed, it can be easily checked that $J(3, 2) \simeq J(3, 1)$.  Thus, the graphs should be $J(2, 1)$ or $J(3, 1)$. For the case $d=2$, the possibilities of the spectrum $T$ are only 
\[ \sigma(T)=\left( -\frac{1}{n-2}, \frac{n-4}{2n-4}, 1 \right)=\left(-1, -\frac{1}{2}, 1\right), \left(-\frac{1}{2}, 0, 1\right). \]
Then it should hold $n=3$ or $4$. If $n=3$, then we have $k=2, 1$ for the cases $d=k, d=n-k$, respectively since $d=2$. However, these pairs do not satisfy the assumption $2k \le n, 2k \ge n$ for the cases $d=k, d=n-k$, respectively. So $n$ should be $4$ in this case. If $n=4$, then we can obtain $k=2$ for both of the cases $d=k,d=n-k$. Thus, the graph should be $J(4, 2)$ in this case.  For the cases $d=3$, and $d=4$, we can obtain 
\[ \sigma(T)=\left( -\frac{1}{n-3}, \frac{n-7}{3n-9}, \frac{2n-9}{3n-9}, 1 \right), \]
\[ \sigma(T)=\left( -\frac{1}{n-4}, \frac{n-10}{4n-16}, \frac{2n-14}{4n-16}, \frac{3n-16}{4n-16}, 1 \right), \]
respectively. However, for every $n$, it follows $\sigma(T) \not\subset \{ \pm 1, \pm 1/2, 0 \}$ for both of the two cases. Thus, there are no periodic Johnson graphs in these cases. So all the periodic Johnson graphs are only $J(2, 1), J(3, 1), J(4, 2)$. 

\end{proof}


\section{Applications to general distance-regular graphs}
\subsection{Necessary condition for periodic distance-regular graphs}
In this section, we consider the periodicity of the Grover walks on the general distance-regular graphs and obtain some necessary conditions. Also, we give another proof of the periodic strongly regular graphs for \cite{Higuchi13}. As is seen in previous sections, the eigenvalue of modified quotient matrix $\tilde{Q}$ should be the real part of a root of unity to induce periodic Grover walks. For a distance-regular graph $G$ with diameter $d$, let $\varphi(x)$ be the characteristic polynomial of its modified quotient matrix $\tilde{Q}$ and $\rho_{i}$ be the coefficient of $x^{i}$ of $\varphi(x)$ for $i \in \{0, \cdots, D\}$, that is, 
\begin{equation}
\varphi(x):=\mathrm{det}(x I_{D}-\tilde{Q})=\sum^{D}_{i=0}\rho_{i}x^{i}, \label{def of phi}
\end{equation} 
where $I_{D}$ is the $D \times D$ identity matrix and we put $D=d+1$ for simplicity. Also, let $\Psi(z)$ be the $2D$-dimensional polynomial applied the Zhukovskij transformation to $\varphi(x)$ and $\alpha_{j}$ be the coefficient of $z^{j}$ of $\Psi(z)$ for $j \in \{0, 1, \cdots, 2D \}$, that is, 
\begin{equation}
\Psi(z):=(2z)^{D} \varphi \left( \frac{z+z^{-1}}{2} \right)=\sum^{2D}_{j=0} \alpha_{j} z^{j}. \label{def of psi}
\end{equation}
Then for the above polynomial, the followings hold. 
\begin{claim}
Let $\alpha_{j}$ be defined as the above. Then it holds that
\begin{description}
\item[(i)]
$\alpha_{2D}=1$, 
\item[(ii)]
$\alpha_{2D-j}=\alpha_{j}$ for every $j \in \{0, 1, \cdots, D \}$. 
\end{description} \label{claim of alpha}
\end{claim}
\begin{proof}
Since $\varphi(x)$ is a monic polynomial, $\Psi(z)$ becomes also a monic polynomial and (i) immediately follows. Applying the Zhukovskij transformation to $\varphi(x)$, we can obtain
\begin{align}
\Psi(z) &=(2z)^{D} \left\{ \sum^{D}_{i=0} \rho_{i} \left(\frac{z+z^{-1}}{2} \right)^{i} \right\}\\
&=\sum^{D}_{i=0} \left\{ 2^{D-i} \rho_{i} \sum^{i}_{r=0} \binom{i}{r}z^{D-2r+i} \right\}. \label{exp of psi}
\end{align}
Considering the coefficient of $z^{2D-j}$ of (\ref{exp of psi}), we have that the term of $i=D-(j-2k)$ with $r=k$ for $k \in \{0, \cdots, \lfloor \frac{j}{2} \rfloor \}$ only contributes to the coefficient of $z^{2D-j}$.  So it follows that 
\begin{equation}
\alpha_{2D-j}=\sum^{\lfloor \frac{j}{2} \rfloor}_{k=0}2^{j-2k} \rho_{D-(j-2k)}  \binom{D-(j-2k)}{k}. \label{alpha2D-j}
\end{equation}
Next, considering the coefficient of $z^{j}$ of (\ref{exp of psi}), we also have that the term of $i=D-(j-2k)$ and $r=D-(j-k)$ for $k \in \{0, \cdots, \lfloor \frac{j}{2} \rfloor \}$ only contributes to $z^{j}$. So we can also obtain 
\begin{equation}
\alpha_{j}=\sum^{\lfloor \frac{j}{2} \rfloor}_{k=0}2^{j-2k} \rho_{D-(j-2k)}  \binom{D-(j-2k)}{D-(j-k)}.
\end{equation}
Since $\binom{D-(j-2k)}{k}=\binom{D-(j-2k)}{D-(j-k)}$, (ii) follows. 
\end{proof}
By Lemma \ref{lem of cyclo}, in order to induce periodic Grover walks, $\Psi(z)$ should be represented by a product of cyclotomic polynomials. So it is required that $\alpha_{j} \in \dZ$ for every $j \in \{0, 1, \cdots, 2D \}$.  
Using this fact, we will prove Theorem \ref{thm of rho}. 

\begin{proof}[Proof of Theorem \ref{thm of rho}]
First, we prove it for every $j \in \{1, \cdots, D-1 \}$ by showing its contraposition.  So we can suppose that there exist $j \in \{1, \cdots, D-1 \}$ such that $2^{j}\rho_{D-j} \notin \dZ$ and $2^{l}\rho_{D-l} \in \dZ$ for every $l \in \{0, 1, \cdots, j-1\}$. By (\ref{alpha2D-j}), we have 
\begin{equation}
 \alpha_{2D-j}=2^{j}\rho_{D-j}+\sum^{\lfloor \frac{j}{2} \rfloor}_{k=1}2^{j-2k} \rho_{D-(j-2k)}  \binom{D-(j-2k)}{k}. 
 \label{alpha2D-j}
 \end{equation}
Then the second term is an integer value and the first term is not from the assumption. Therefore, $\alpha_{2D-j} \notin \dZ$, which implies that $\Psi(z)$ cannot be represented by a product of cyclotomic polynomials then $G$ is not periodic by Lemma \ref{lem of cyclo}.  

Next, we prove why we can omit the case for $j=D$. Since $1$ is always an eigenvalue of $T$, $\varphi(x)$ always has a factor $(x-1)$. Considering the Zhukovskij transformation, it follows that $\Psi(z)$ always has a factor $(z-1)$. So we have $\Psi(1)=0$, that is, 
\begin{equation}
\Psi(1)=\sum^{2D}_{j=0}\alpha_{j}=0.
\end{equation}
Then using (ii) of Claim \ref{claim of alpha}, we can also obtain 
\begin{equation}
\Psi(1)=2\sum^{D-1}_{j=0} \alpha_{2D-j}+\alpha_{D}=0. \label{alphaD}
\end{equation}
If it holds $2^{j}\rho_{D-j} \in \dZ$ for every $j \in \{1, \cdots, D-1 \}$, then the first term is an integer value from (\ref{alpha2D-j}). By the above equation, we have $\alpha_{D} \in \dZ$, which implies that $2^{D}\rho_{0} \in \dZ$ from (\ref{alpha2D-j}). Therefore, if the statement holds for every $j \in \{1, \cdots, D-1 \}$, then the case for $j=D$ should be an integer value.
\end{proof}

Indeed, Theorem \ref{thm of rho} can be extended to a general graph since we just used the coefficients of its transition matrix. So we can show the following by the same procedure. 
\begin{rem}
For a graph $G$ with $|G|=n$ and its transition matrix $T$, let 
\[ \varphi(x)=\mathrm{det}(xI_{n}-T):=\sum^{n}_{i=0}\rho_{i}x^{i}. \]
If $G$ is periodic, then it should hold
\[ 2^{j}\rho_{n-j} \in \dZ \]
for every $j \in \{ 0, 1, \cdots, n \}$. 
\end{rem}

\subsection{Application of Theorem \ref{thm of rho}: Another proof of Proposition \ref{prop srg}}
From the above arguments, analyzing the coefficients of $\varphi(x)$ is one of the method to consider its periodicity. Let us denote the matrix $x I_{D}-\tilde{Q}$ by $Y$. From the definition of the determinant, we have
\begin{equation}
 \mathrm{det}Y=\sum_{\sigma \in S_{D}} \mathrm{sgn}(\sigma) Y_{0, \sigma(0)}Y_{1, \sigma(1)}\cdots Y_{d, \sigma(d)}, \label{def of det}
\end{equation}
where $S_{D}$ is the symmetric group on $\{0, 1, \cdots, d \}$ and $\sigma$ corresponds to a permutation on $\{0, 1, \cdots, d \}$. So we can calculate the coefficients of $\varphi(x)$ by (\ref{def of det}).

Here, we consider the periodic distance-regular graph with $d=2$, which is nothing but a strongly regular graph. However, the periodic strongly regular graphs have been already known. These are only $C_{5}, K_{k, k}, K_{\lambda, \lambda, \lambda}$ for $k, \lambda \in \dN$ as is seen in Proposition \ref{prop srg}. So we provide another proof for periodic strongly regular graphs by using Theorem \ref{thm of rho}. 
\begin{thm}
The periodic distance regular graphs with $d=2$ are only $C_{5}, K_{k, k}, K_{m, m, m}$ for $k, m \in \dN$.
\end{thm}
\begin{proof}[Another proof]
From Theorem \ref{thm of rho},  it should hold that
\begin{equation}
2\rho_{2}, 4\rho_{1} \in \dZ \label{con of rho}
\end{equation}
due to the periodicity. First, we consider $\rho_{2}$. The permutation on (\ref{def of det}) which contributes to the term $x^{2}$ of $\varphi(x)$ is only the identity permutation, whose sign is $+1$ . Then it can easily follow that $\rho_{2}$ is constructed by the summation of all the products of two variables $x$ and one constant value $-a_{i}/k$ for $i \in \{0, 1, 2\}$, that is, 
\[ \rho_{2}=-\left(\frac{a_{0}}{k}+\frac{a_{1}}{k}+\frac{a_{2}}{k} \right). \]
By $a_{0}=0$ and (\ref{con of rho}), we have
\begin{equation}
-2\left(\frac{a_{1}}{k}+\frac{a_{2}}{k} \right) \in \dZ. \label{rho2}
\end{equation}
Next, we consider $\rho_{1}$. The permutations on (\ref{def of det}) which contribute to the term $x$ of $\varphi(x)$ are only the identity permutation and the transposition on $\{i, j \}$ with $\tilde{Q}_{i, j}\ne 0$, that is, $j=i+1$ for $i \in \{0, 1 \}$. Then their signs are $+1$, $-1$, respectively. For the identity permutation, we can choose one variable $x$ and two constant values $-a_{i}/k, -a_{j}/k$ with $i, j \in \{0, 1, 2 \}$, $i \ne j$ and for a transposition, we can choose one variable $x$ and two constant values $-b_{i}/k, -c_{i+1}/k$ with $i \in \{0, 1 \}$. Thus, it holds 
\[
\rho_{1}= \left(\frac{a_{0}a_{1}}{k^{2}}+\frac{a_{0}a_{1}}{k^{2}}+\frac{a_{1}a_{2}}{k^{2}} \right)-\left(\frac{b_{0}c_{1}}{k^{2}} +\frac{b_{1}c_{2}}{k^{2}} \right).
\] 
By $a_{0}=0$ and (\ref{con of rho}), we can also obtain
\begin{equation}
4\left\{ \frac{a_{1}a_{2}}{k^{2}} -\left(\frac{b_{0}c_{1}}{k^{2}} +\frac{b_{1}c_{2}}{k^{2}} \right)  \right\} \in \dZ. \label{rho1}
\end{equation}
Using (\ref{rho2}), (\ref{rho1}), we characterize these parameters. 

First, we treat the case $a_{1}+a_{2}=0$, which implies that $a_{1}=a_{2}=0$ and it obviously satisfies the condition (\ref{rho2}).  From the definition of distance-regular graphs, the values $b_{0}$ and $c_{1}$ are always $k, 1$, respectively. Combining them and (\ref{c+a+b}), it follows $b_{1}=k-1$ and $c_{2}=k$. So it has been already obtained that $a_{1}=a_{2}=0$, and $b_{0}=k, b_{1}=k-1$, and $c_{1}=1, c_{2}=k$. Assigning them to (\ref{rho1}), it becomes $-4 \in \dZ$ for any $k \in \dN$. Therefore, the quotient matrix $Q$ becomes
\begin{equation}
Q=\left(
	\begin{array}{ccc}
	0 & k & 0 \\
	1 & 0 & k-1\\
	0 & k & 0
	 \end{array}
	\right).  \label{Kkk}
\end{equation}
Then we can calculate all the coefficients of $\Psi(z)$ by (\ref{alpha2D-j}), (\ref{alphaD}) as follows: $\alpha_{6}=\alpha_{0}=1, \alpha_{5}=\alpha_{1}=0, \alpha_{4}=\alpha_{2}=-1, \alpha_{3}=0$ from the above $Q$ and we have
 \[ \Psi(z)=z^{6}-z^{4}-z^{2}+1=(z^{2}-1)(z^{4}-1), \]
whose roots satisfy $z^{4}=1$. Therefore, in this case, the distance-regular graph is periodic for any $k$. Indeed, we can obtain that the distance-regular graph which achieves the quotient matrix (\ref{Kkk}) is uniquely determined as $K_{k, k}$. (See Appendix)

From now on, we will suppose $a_{1}+a_{2}\ne 0$. If $k$ is odd, then it should hold 
\[ a_{1}+a_{2} \in k\dZ \]
by (\ref{rho2}). Since $0 \le a_{i} < k$, we have $0 < a_{1}+a_{2} <2k$. So we can suppose $a_{1}+a_{2}=k$. From this assumption and (\ref{c+a+b}), it follows $c_{2}=a_{1}, b_{1}=a_{2}-1$. Assigning them and $b_{0}=k, c_{1}=1$ to (\ref{rho1}), we have
\begin{align*}
4\left\{ \frac{a_{1}a_{2}}{k^{2}} -\left(\frac{b_{0}c_{1}}{k^{2}} +\frac{b_{1}c_{2}}{k^{2}} \right)  \right\} &= 4\left\{ \frac{a_{1}a_{2}}{k^{2}}-\left(\frac{k+(a_{2}-1)a_{1}}{k^{2}}\right) \right\} \\
&= 4\left( \frac{a_{1}-k}{k^{2}} \right). 
\end{align*} 
Since $k$ is odd, it should hold $a_{1}-k \in k^{2}\dZ$. However, there is the restriction $0 \le a_{1} < k$. So we have $-k \le a_{1}-k < 0$. Hence, $a_{1}-k \notin k^{2}\dZ$. So there is no periodic distance-regular graph under this condition. 

If $k$ is even, that is, there exists $m \in \dN$ such that $k=2m$, then it should hold
\[ a_{1}+a_{2} \in m\dZ \]
by (\ref{rho2}). If $m=1$, then $k$, the valency of $G$, is $2$. Indeed, the connected $2$-regular graph with $d=2$ is only $C_{5}$ and it is a well-known periodic graph \cite{Yoshie}. So we can assume $m>1$. Since $0<a_{1}+a_{2}< 2k=4m$, the candidates are only $a_{1}+a_{2}=m, 2m, 3m$. If $a_{1}+a_{2}=m$, then we have $b_{1}=2m-a_{1}-1, c_{2}=m-a_{1}$ by the assumption and (\ref{c+a+b}). Assigning them and $b_{0}=k=2m, c_{1}=1$ to (\ref{rho1}), it is required that the following
\begin{align*}
4\left\{ \frac{a_{1}a_{2}}{k^{2}} -\left(\frac{b_{0}c_{1}}{k^{2}} +\frac{b_{1}c_{2}}{k^{2}} \right)  \right\}  &= -\left\{ \frac{a_{1}a_{2}}{m^{2}}-\left( \frac{b_{0}c_{1}+b_{1}c_{2}}{m^{2}} \right) \right\} \\
&= -\left\{ \frac{a_{1}(m-a_{1})}{m^{2}}-\left( \frac{2m+(2m-a_{1}-1)(m-a_{1})}{m^{2}} \right) \right\}\\
&= -\frac{1}{m^{2}}\left( -2m^{2}+a_{1}-m \right)\\
&= 2-\frac{1}{m^{2}}(a_{1}-m)
\end{align*}
should be an integer value, that is, $a_{1}-m \in m^{2}\dZ$. Due to the restriction $0 < a_{1} < k=2m$ and $m > 1$, we have 
\[ -m^{2}+m < -m <a_{1}-m < m < m^{2}+m. \]
Hence, it should hold $a_{1}-m=0$, that is, $a_{1}=m$. Then it immediately follows $a_{2}=0, b_{1}=m-1, c_{2}=2m$ and its quotient matrix $Q$ becomes 
\begin{equation}
Q=\left(
	\begin{array}{ccc}
	0 & 2m & 0 \\
	1 & m & m-1\\
	0 & 2m & 0
	 \end{array}
	\right).  \label{Kmmm}
\end{equation}
Then we can calculate all the coefficients of $\Psi(z)$ by (\ref{alpha2D-j}), (\ref{alphaD}) as follows: $\alpha_{6}=\alpha_{0}=1, \alpha_{5}=\alpha_{1}=-1, \alpha_{4}=\alpha_{2}=1, \alpha_{3}=-2$ by the above $Q$ and we have
 \[ \Psi(z)=z^{6}-z^{5}+z^{4}-2z^{3}+z^{2}-z+1=(z-1)(z^{2}+1)(z^{3}-1), \]
whose roots satisfy $z^{12}=1$.
So we have that in this case, the distance-regular graph whose quotient matrix can be represented by (\ref{Kmmm}) is only periodic. Indeed, such distance-regular graph can be uniquely determined as $K_{m, m, m}$ (See Appendix). If $a_{1}+a_{2}=2m$ or $3m$, then by similar method, we can also obtain the conditions 
\[ \frac{a_{1}-2m}{m^{2}}, \frac{a_{1}-3m}{m^{2}} \in \dZ, \]
respectively. However these values cannot be integer values due to the restriction $0<a_{1}<k=2m$. So there is no periodic distance-regular graph in these cases. Therefore, it can be also proven that the periodic distance-regular graphs with $d=2$ are only $C_{5}, K_{k, k}, K_{m, m, m}$ for $k, m \in \dN$. 
\end{proof}
 
\section{Summary and discussions}
In this paper, we characterized all of periodic Hamming graphs and Johnson graphs. The key tool for such characterizations is Theorem \ref{thm of diam}.  It is useful for every graphs with $\sigma(T) \subset \dQ$. So analyzing the periodicity of the Grover walks on such graphs, we can control the diameter of the graphs. Therefore, characterizing such graphs is also our future work. In the latter part of this paper, we obtained a necessary condition (Theorem \ref{thm of rho}) for general periodic distance-regular graphs. By using it, we could classify all the periodic distance-regular graphs with $d=2$, that is, we provided another proof for periodic strongly regular graphs. For $d \ge 3$, the computation becomes harder since the number of parameters increases. So we should impose another condition except for Theorem $\ref{thm of rho}$ to characterize all the periodic distance-regular graphs. To find such condition is also our future work and we would like to characterize all the periodic distance-regular graphs by mixing them. Also, Theorem \ref{thm of rho} can be extended to every graphs because it requires only the coefficients of characteristic polynomial of its transition matrix. Indeed, these coefficients related to cycles or matchings on the graph. Then we can control the shape of graphs by only supposing the periodicity of the Grover walk. So we have to classify another property and determine the shape of the graphs by it. 

\section*{Acknowledgments}
The author thanks to Etsuo Segawa and Yusuke Higuchi for supporting our researches and giving us fruitful comments.

\section*{Appendix}
In this part, we prove that distance-regular graphs with quotient matrices (\ref{Kkk}), (\ref{Kmmm}) are uniquely determined as $K_{k, k}$, and $K_{m, m, m}$, respectively. 

First, we consider (\ref{Kkk}):
\begin{equation*}
Q=\left(
	\begin{array}{ccc}
	0 & k & 0 \\
	1 & 0 & k-1\\
	0 & k & 0
	 \end{array}
	\right). 
\end{equation*}
Let $G$ be a distance-regular graph with quotient matrix (\ref{Kkk}). Then $G$ is a $k$-regular graph. We fix $x \in V(G)$. Then we have  $|\Gamma_{0}(x)|=1, |\Gamma_{1}(x)|=k, |\Gamma_{2}(x)|=k-1$ by the relation 
\begin{equation}
b_{j} |\Gamma_{j}(x)|=c_{j+1} |\Gamma_{j+1}(x)| \label{bc}
\end{equation}

for $j \in \{0, 1, 2\}$ \cite{GM} and $|\Gamma_{0}(x)|=1$, where $|\Gamma_{j}(x)|$ is defined in section 2. So any vertices in $\Gamma_{1}(x)$ are adjacent to all the vertices in $\Gamma_{2}(x)$ and vice versa since $b_{1}=k-1, c_{2}=k$. Put 
\begin{align*}
 V_{1}=&\Gamma_{0}(x) \cup \Gamma_{2}(x),\\
 V_{2}=&\Gamma_{1}(x). 
\end{align*}
Then $|V_{1}|=|V_{2}|=k$, and these two partitions become complete bi-partitions. Therefore, $G$ can be determined as $K_{k, k}$. 

Next, we consider (\ref{Kmmm}):
\begin{equation*}
Q=\left(
	\begin{array}{ccc}
	0 & 2m & 0 \\
	1 & m & m-1\\
	0 & 2m & 0
	 \end{array}
	\right). 
\end{equation*}
Let $G$ be a distance-regular graph with quotient matrix (\ref{Kmmm}). Then $G$ is a $2m$-regular graph. Also, we fix $x \in V(G)$ and it follows $|\Gamma_{0}(x)|=1, |\Gamma_{1}(x)|=2m, |\Gamma_{2}(x)|=m-1$ from (\ref{bc}).  For a vertex $u \in \Gamma_{1}(x)$, put 
\begin{align*}
V_{1}=&\Gamma_{0}(x) \cup \Gamma_{2}(x),\\
V_{2}=&N(u)\cap \Gamma_{1}(x),\\
V_{3}=&\Gamma_{1}(x) \backslash V_{2}. 
\end{align*}
Then we have $|V_{1}|=|V_{2}|=|V_{3}|=m$ since $|\Gamma_{2}(x)|=m-1$, $a_{1}=m$ and $u \in V_{3}$. We show these become complete tri-partitions. From the definition of $\Gamma_{1}(x)$,  $x$ is adjacent to all the vertices in both of $V_{2}$ and $V_{3}$. In addition, any vertices in $\Gamma_{2}(x)$ are adjacent to all the vertices in $\Gamma_{1}(x)$ since $c_{2}=|\Gamma_{1}(x)|=2m$. So we can obtain that any vertices in $V_{1}$ are adjacent all the vertices in $V_{2} \cup V_{3}$. Then no two vertices in $V_{1}$ are adjacent each other and we will show it for both of $V_{2}, V_{3}$. Suppose there exists $u_{1}, u_{2} \in V_{2}$ such that $u_{1} \sim u_{2}$. Since $|V_{3}|=m$ and $m$ edges from $u_{1}$ extend to $V_{1}$, there exists a vertex $w(\ne u) \in V_{3}$ such that $u_{1} \not\sim w$. Also, it follows $u \not\sim w$ since $m$ edges extend from $u$ to $V_{1}$ and the remnant $m$ edges extend from $u$ to $V_{2}$.  Now we retake $\Gamma_{0}(w), \Gamma_{1}(w), \Gamma_{2}(w)$. Since both of $u, u_{1}$ are not adjacent to $w$, it follows $u, u_{1} \in \Gamma_{2}(w)$. Moreover, it follows $u \sim u_{1}$ from the definition of $V_{2}$, which contradicts $a_{2}=0$. So we have that any vertices in $V_{2}$ are not adjacent to another vertex in $V_{2}$ and it also follows for $V_{3}$ since $|V_{3}|=m$. Hence, these three become complete tri-partitions and $G$ can be uniquely determined as $K_{m, m, m}$.

\end{document}